\documentclass[12pt, reqno]{amsart}
\usepackage{amsmath, amsthm, amscd, amsfonts, amssymb, graphicx, color}
\usepackage[bookmarksnumbered, colorlinks, plainpages]{hyperref}
\hypersetup{colorlinks=true,linkcolor=red, anchorcolor=green, citecolor=cyan, urlcolor=red, filecolor=magenta, pdftoolbar=true}

\newtheorem{theorem}{Theorem}[section]
\newtheorem{lemma}[theorem]{Lemma}
\newtheorem{proposition}[theorem]{Proposition}
\newtheorem{corollary}[theorem]{Corollary}
\theoremstyle{definition}

\theoremstyle{remark}

\numberwithin{equation}{section}

\begin{document}

\setcounter{page}{1}

\title[Hardy and Rellich inequalities]{Hardy and Rellich inequalities for anisotropic $p$-sub-Laplacians}

\author[M. Ruzhansky, B. Sabitbek, \MakeLowercase{and} D. Suragan]{M. Ruzhansky,$^1$ B. Sabitbek,$^{2^*}$ \MakeLowercase{and}  D. Suragan$^{3}$}

\address{$^{1}$Department of Mathematics: Analysis, Logic and Discrete Mathematics, Ghent University, Krijgslaan 281, Building S8 
	B 9000, Ghent, Belgium, and School of Mathematical Sciences, Queen Mary University of London, United Kingdom}
\email{\textcolor[rgb]{0.00,0.00,0.84}{Michael.Ruzhansky@ugent.be}}

\address{$^{2}$Institute of Mathematics and Mathematical Modeling, 125 Pushkin str, 050010, Almaty, Kazakhstan
	and Al-Farabi Kazakh National University, Almaty, Kazakhstan}
\email{\textcolor[rgb]{0.00,0.00,0.84}{b.sabitbek@math.kz}}

\address{$^{3}$Department of Mathematics, Nazarbayev University,
	53 Kabanbay Batyr Ave, 010000, Astana, Kazakhstan}
\email{\textcolor[rgb]{0.00,0.00,0.84}{durvudkhan.suragan@nu.edu.kz}}

\let\thefootnote\relax\footnote{Copyright 2018 by the Tusi Mathematical Research Group.}



\date{Received: xxxxxx; Revised: yyyyyy; Accepted: zzzzzz.
\newline \indent $^{*}$Corresponding author}

\begin{abstract}
In this paper we establish the subelliptic Picone type identities. As consequences, we obtain Hardy and Rellich type inequalities for anisotropic $p$-sub-Laplacians which are operators of the form
\begin{equation*}
	\mathcal{L}_{p}f:= \sum_{i=1}^{N} X_i\left( |X_i f|^{p_i-2} X_i f  \right) ,\quad 1<p_i<\infty,
\end{equation*}
where $X_i$, $i=1,\ldots, N$, are the generators of the first stratum of a stratified (Lie) group.
Moreover, analogues of Hardy type inequalities with multiple singularities and many-particle Hardy type inequalities are obtained on stratified groups.
\end{abstract} \maketitle

\section{\textbf{Introduction}}
	\subsection{Historical background}
Recall the anisotropic Laplacian on $\mathbb R^N$ for $p_i>1$ where $ i=1,\dots,N$ (see \cite{Feng-Cui}), defined by
\begin{equation}\label{eq1}
\sum_{i=1}^{N} \frac{\partial }{\partial x_i} \left(\left|\frac{\partial u}{\partial x_i}\right|^{p_i-2} \frac{\partial u}{\partial x_i}\right),
\end{equation}
which has recently attracted considerable attention. For instance, by taking $p_i=2$ or $p_i=p=const$ (see \cite{BK}) in \eqref{eq1}  we get the Laplacian and the pseudo-$p$-Laplacian, respectively. The anisotropic Laplacian has the theoretical importance not only in mathematics, but also many practical applications in the natural sciences. There are several examples: it reflects anisotropic physical properties of some reinforced materials (Lions \cite{Lions} and Tang \cite{Tang}), as well as explains the dynamics of fluids in the anisotropic media when the conductivities of the media are different in each direction \cite{ADSh} and \cite{Bear}. It has also applications in the image processing \cite{Weickert}.

The main purpose of this paper is to obtain the subelliptic Hardy and Rellich type inequalities for the anisotropic $p$-sub-Laplacian on stratified groups. First, we derive the subelliptic Picone type identities on stratified groups. As consequences, Hardy and Rellich type inequalities for anisotropic sub-Laplacians are presented. These results are given in Sections \ref{sec1-1} and \ref{sec2}. In Section \ref{sec3} and \ref{sec4}, we present analogues of Hardy type inequalities with multiple singularities and many-particle Hardy type inequalities on stratified groups, respectively. These inequalities are obtained in their horizontal form in the spirit of \cite{RS17a} (see also \cite{RS18}). In Section \ref{sec5}, Hardy type inequalities with exponential weights are shown.
\subsection{Preliminaries}
Let $\mathbb{G}=(\mathbb{R}^d,\circ,\delta_{\lambda})$ be a stratified Lie group (or a homogeneous Carnot group), with dilation structure $\delta_{\lambda}$ and Jacobian generators $X_{1},\ldots,X_{N}$, so that $N$ is the dimension of the first stratum of $\mathbb{G}$. We denote by $Q$ the homogeneous dimension of $\mathbb{G}$.  We refer to \cite{PR01, PR08} and to the recent book \cite{FR} for extensive discussions of stratified Lie groups and their properties.

The sub-Laplacian on $\mathbb{G}$ is given by
\begin{equation}
\mathcal{L}=\sum_{k=1}^{N}X_{k}^{2}.
\end{equation}
We also recall that the standard Lebesgue measure $dx$ on $\mathbb R^{n}$ is the Haar measure for $\mathbb{G}$ (see, e.g. \cite{FR}).
The left invariant vector field $X_{j}$ has an explicit form and satisfies the divergence theorem,
see e.g. \cite{FR} for the derivation of the exact formula: more precisely, we can formulate 
\begin{equation}
X_{k}=\frac{\partial}{\partial x'_{k}}+
\sum_{l=2}^{r}\sum_{m=1}^{N_{l}}a_{k,m}^{(l)}(x',...,x^{(l-1)})
\frac{\partial}{\partial x_{m}^{(l)}},
\end{equation}
with $x=(x',x^{(2)},\ldots,x^{(r)})$, where $r$ is the step of $\mathbb{G}$ and
$x^{(l)}=(x^{(l)}_1,\ldots,x^{(l)}_{N_l})$ are the variables in the $l^{th}$ stratum.
The horizontal gradient is given by
$$\nabla_{\mathbb{G}}:=(X_{1},\ldots, X_{N}),$$
and the horizontal divergence is defined by
$${\rm div}_{\mathbb{G}} v:=\nabla_{\mathbb{G}}\cdot v.$$
The horizontal anisotropic $p$-sub-Laplacian is defined by
\begin{equation}
\mathcal{L}_{p}f:= \sum_{i=1}^{N} X_i\left( |X_i f|^{p_i-2} X_i f  \right) ,\quad 1<p_i<\infty,
\end{equation}
and we use the notation 
$$|x'|=\sqrt{x'^{2}_{1}+\ldots+x'^{2}_{N}}$$ for the Euclidean norm on $\mathbb{R}^{N}.$

\section{\textbf{Subelliptic anisotropic Hardy type inequality}}
\label{sec1-1}
First, we obtain the subelliptic Picone type identity on a stratified group $\mathbb{G}$. Here we follow the method of Allegretto and Huang \cite{AH} for the (Euclidean) $p$-Laplacian  (see also \cite{PHY}).
\begin{lemma}\label{Picone}
	Let $\Omega \subset \mathbb{G}$ be an open set, where $\mathbb{G}$ is a  stratified group with $N$ being the dimension of the first stratum. Let $u,v$ be differentiable a.e. in $\Omega$, $v>0$ a.e. in $\Omega$ and $u\geq 0$, and denote
	\begin{equation*}
	R(u,v) := \sum_{i=1}^{N} \left| X_i u\right|^{p_i} - \sum_{i=1}^{N} X_i \left(\frac{u^{p_i}}{v^{p_i-1}}\right)\left|X_i v \right|^{p_i-2}X_iv,
	\end{equation*}
	\begin{align*}
	L(u,v) := \sum_{i=1}^{N} \left| X_i u \right|^{p_i} &- \sum_{i=1}^{N}p_i\frac{u^{p_i-1}}{v^{p_i-1}} \left|X_i v\right|^{p_i-2} X_ivX_i u \nonumber\\
	&+ \sum_{i=1}^{N} (p_i-1)\frac{u^{p_i}}{v^{p_i}} \left|X_i v\right|^{p_i},
	\end{align*}
	where $p_i>1$, $i=1,\ldots,N$. Then
	\begin{equation}\label{1.3}
	L(u,v)=R(u,v) \geq 0.
	\end{equation}
	In addition, let $\Omega$ be connected, then we have $L(u,v)=0$ a.e. in $\Omega$ if and only if $u=cv$ a.e. in $\Omega$ with a positive constant $c$.
\end{lemma}
Note that the Euclidean case of this lemma was obtained by Feng and Cui  \cite{Feng-Cui}.
\begin{proof}[Proof of Lemma \ref{Picone}]
	A direct computation gives
	\begin{align*}
	R(u,v) &= \sum_{i=1}^{N} \left|X_i u\right|^{p_i} - \sum_{i=1}^{N} X_i \left( \frac{u^{p_i}}{v^{p_i-1}}\right) |X_i v|^{p_i-2}X_i v \\
	&= \sum_{i=1}^{N} \left|X_i u\right|^{p_i} - \sum_{i=1}^{N}p_i\frac{u^{p_i-1}}{v^{p_i-1}}|X_iv|^{p_i-2}X_i v X_iu +\sum_{i=1}^{N} (p_i-1)\frac{u^{p_i}}{v^{p_i}}|X_iv|^{p_i}\\
	& = L(u,v).
	\end{align*}
	This proves the equality in \eqref{1.3}. Now we rewrite $L(u,v)$ to see $L(u,v)\geq 0$, that is,
	\begin{align*}
	L(u,v) = & \sum_{i=1}^{N} |X_i u|^{p_i} - \sum_{i=1}^{N} p_i \frac{u^{p_i-1}}{v^{p_i-1}} |X_i v|^{p_i-1}|X_i u| + \sum_{i=1}^{N}(p_i-1)\frac{u^{p_i}}{v^{p_i}}|X_i v|^{p_i} \\
	+&\sum_{i=1}^{N}p_i \frac{u^{p_i-1}}{v^{p_i-1}} |X_i v|^{p_i-2} \left(|X_i v||X_i u| -X_i v X_iu  \right) \\
	=& S_1 +S_2,
	\end{align*}
	where we denote
	\begin{align*}
	S_1:=& \sum_{i=1}^{N} p_i \left[ \frac{1}{p_i} |X_i u|^{p_i} + \frac{p_i-1}{p_i}\left(\left(\frac{u}{v}|X_i v|\right)^{p_i-1}\right)^{\frac{p_i}{p_i-1}}\right] \\
	-& \sum_{i=1}^{N}p_i \frac{u^{p_i-1}}{v^{p_i-1}}|X_i v|^{p_i-1}|X_i u|,
	\end{align*}
	and
	\begin{equation*}
	S_2:= \sum_{i=1}^{N}p_i \frac{u^{p_i-1}}{v^{p_i-1}}|X_i v|^{p_i-2}\left(|X_i v||X_i u| -X_i v X_iu\right).
	\end{equation*}
	We can see that $S_2\geq 0$ due to $|X_iv||X_i u| \geq X_i v X_iu$. To check $S_1 \geq 0$ we need to use Young's inequality for $a\geq 0$ and $b \geq 0$
	\begin{equation}\label{2.2}
	ab \leq \frac{a^{p_i}}{p_i} + \frac{b^{q_i}}{q_i},
	\end{equation}
	where $p_i>1, q_i>1$ and $\frac{1}{p_i}+\frac{1}{q_i}=1$ for $i=1,\ldots,N$. It holds if and only if $a^{p_i}=b^{q_i}$, i.e. if $a = b^{\frac{1}{p_i-1}}$. Let us take $a= |X_i u|$ and $b =\left(\frac{u}{v}|X_iv|\right)^{p_i-1}$ in \eqref{2.2} to get
	\begin{equation}
	p_i |X_i u|\left(\frac{u}{v}|X_i v|\right)^{p_i-1} \leq p_i \left[\frac{1}{p_i} |X_i u|^{p_i} + \frac{p_i-1}{p_i}\left(\left(\frac{u}{v}|X_i v|\right)^{p_i-1}\right)^{\frac{p_i}{p_i-1}}\right].
	\end{equation}
	From this we see that $S_1 \geq 0$ which proves that $L(u,v)=S_1+S_2 \geq 0$.
	It is easy to see that $u=cv$ implies $R(u,v)=0$ since in the case $R(cv,v)$ both terms cancel each other out. Now let us prove that $L(u,v)=0$ implies $u=cv$. Due to $u(x)\geq 0$ and $L(u,v)(x_0)=0, \quad x_0\in \Omega,$ we consider the two cases $u(x_0)>0$ and $u(x_0) =0.$ 
	\begin{itemize}
		\item[(1)] For the case $u(x_0)>0$ we conclude from $L(u,v)(x_0)=0$ that $S_1=0$ and $S_2=0$. Then $S_1=0$ implies
		\begin{equation}\label{2.4}
		|X_i u| = \frac{u}{v} |X_i v|, \quad i=1,\ldots,N,
		\end{equation}
		and $S_2=0$ implies
		\begin{equation}\label{eee}
		|X_i v| |X_i u| - X_i v X_i u = 0, \quad i=1,\ldots,N.
		\end{equation}
		The combination of \eqref{2.4} and \eqref{eee} gives
		\begin{equation}
		\frac{X_i u}{X_i v} = \frac{u }{v}=c, \quad \text{with} \quad c\neq 0, \quad i=1,\ldots,N .
		\end{equation}
		
		\item[(2)]  Let us denote $\Omega^*:=\{x \in \Omega | u(x)=0 \}$. If $\Omega^* \neq \Omega$, then suppose that $x_0 \in \partial \Omega^*$. Then there exists a sequence $x_k \notin \Omega^*$ such that $x_k \rightarrow x_0$. In particular, $u(x_k) \neq 0$, and hence by the case 1 we have $u(x_k)=cv(x_k)$. Passing to the limit we get $u(x_0) = c v(x_0)$. Since $u(x_0)= 0, \quad v(x_0)\neq 0$, we get that $c=0$. But then by the case 1 again, since $u =cv$ and $u\neq 0$ in $\Omega \backslash \Omega^*$, it is impossible that $c=0$. This contradiction implies that $\Omega^* = \Omega$.   
	\end{itemize}
	This completes the proof of Lemma \ref{Picone}.	
\end{proof}
As a consequence of the subelliptic Picone type identity, we present the Hardy type inequality for the anisotropic sub-Laplacian on $\mathbb{G}$. We recall that for $x\in \mathbb{G}$ we write $x=(x',x'')$, with $x'$ being in the first stratum of $\mathbb{G}$.
\begin{theorem}\label{Hardy}
	Let $\Omega \subset \mathbb{G} \backslash \{x'=0\}$ be an open set, where $\mathbb{G}$ is a stratified group  with $N$ being the dimension of the first stratum.
	Then we have
	\begin{equation}\label{1.8}
	\sum_{i=1}^{N} \int_{\Omega} |X_i u|^{p_i} dx \geq \sum_{i=1}^{N}\left(\frac{p_i-1}{p_i}\right)^{p_i} \int_{\Omega} \frac{|u|^{p_i}}{|x'_i|^{p_i}}dx,
	\end{equation}
	for all $u \in C^1(\Omega)$ and where $1<p_i<N$ for $i=1,\ldots,N.$
\end{theorem}
Before we start the proof of Theorem \ref{Hardy}, let us establish the following  Lemma \ref{lem1}.
\begin{lemma}\label{lem1}
	Let $\Omega \subset \mathbb{G}$ be an open set, where $\mathbb{G}$ is a stratified group  with $N$ being the dimension of the first stratum. Let constants $K_i>0$ and functions $H_i(x)$ with $i=1,\ldots,N$, be such that for an a.e. differentiable function $v$, such that $v >0$ a.e. in $\Omega$, we have
	\begin{equation}\label{3.1}
	- X_i (|X_i v|^{p_i-2}X_i v)\geq K_i H_i(x) v^{p_i-1},\quad i=1,\ldots,N.
	\end{equation}
	Then, for all nonnegative functions $u \in C^1(\Omega)$ we have
	\begin{equation}
	\sum_{i=1}^{N}\int_{\Omega} |X_i u|^{p_i} dx \geq 	\sum_{i=1}^{N} K_i \int_{\Omega} H_i(x) u^{p_i} dx.
	\end{equation}
\end{lemma}
\begin{proof}[Proof of Lemma \ref{lem1}.]
	In view of \eqref{1.3} and \eqref{3.1} we have
	\begin{align*}
	0 \leq \int_{\Omega} L(u,v) dx &= \int_{\Omega} R(u,v) dx \\
	& = \sum_{i=1}^{N} \int_{\Omega} |X_i u|^{p_i} dx - 	\sum_{i=1}^{N} \int_{\Omega}X_i \left(\frac{u^{p_i}}{v^{p_i-1}}\right)|X_i v|^{p_i-2}X_iv dx\\
	& = \sum_{i=1}^{N} \int_{\Omega} |X_i u|^{p_i} dx + 	\sum_{i=1}^{N} \int_{\Omega} \frac{u^{p_i}}{v^{p_i-1}} X_i\left(|X_i v|^{p_i-2}X_iv\right) dx\\
	& \leq \sum_{i=1}^{N} \int_{\Omega} |X_i u|^{p_i} dx -  \sum_{i=1}^{N} K_i\int_{\Omega} H_i(x) u^{p_i} dx.
	\end{align*}
	This completes the proof of Lemma \ref{lem1}.
\end{proof}
\begin{proof}[Proof of Theorem \ref{Hardy}.]
	Before using Lemma \ref{lem1}, we shall introduce the auxiliary function
	\begin{equation}
	v := \prod_{j=1}^{N}|x'_j|^{\alpha_j} = |x'_i|^{\alpha_i} V_i,
	\end{equation}
	where $V_i = \prod_{j=1,j\neq i}^{N}|x_j'|^{\alpha_j}$ and $\alpha_j = \frac{p_j-1}{p_j}$. Then we have
	\begin{align*}
	X_i v &= \alpha_i V_i |x_i'|^{\alpha_i -2} x'_i, \\
	|X_i v|^{p_i-2} &= \alpha_i^{p_i-2} V_i^{p_i-2} |x_i'|^{\alpha_ip_i - 2\alpha_i - p_i+2},	\\
	|X_i v|^{p_i-2}X_i v &= \alpha_i^{p_i-1} V_i^{p_i -1}|x_i'|^{\alpha_i p_i - \alpha_i -p_i} x_i'.	
	\end{align*}
	
	Consequently, we also have
	\begin{equation}
	- X_i (|X_i v|^{p_i-2}X_i v) = \left(\frac{p_i-1}{p_i}\right)^{p_i} \frac{v^{p_i-1}}{|x_i'|^{p_i}}.
	\end{equation}
	To complete the proof of Theorem \ref{Hardy}, we choose $K_i = \left(\frac{p_i-1}{p_i}\right)^{p_i}$ and $H_i(x) = \frac{1}{|x_i'|^{p_i}}$, and use Lemma \ref{lem1}.
\end{proof}

\section{\textbf{Subelliptic anisotropic Rellich type inequality}}\label{sec2}

We now present the (second order) subelliptic Picone type identity. As a consequence, we  obtain the Rellich type inequality for the anisotropic sub-Laplacian on $\mathbb{G}$.
\begin{lemma}\label{Picone1}
	Let $\Omega \subset \mathbb{G}$ be an open set, where $\mathbb{G}$ is a stratified group  with $N$ being the dimension of the first stratum. Let $u,v$ be twice differentiable a.e. in $\Omega$ and satisfying the following conditions: $u\geq0$, $v>0$, $X_i^2v<0$ a.e. in $\Omega$ for $p_i>1$, $i=1,\ldots,N$. Then we have
	\begin{equation}\label{L=R}
	L_1(u,v) = R_1(u,v) \geq 0,
	\end{equation}
	where
	\begin{equation*}
	R_1(u,v) := \sum_{i=1}^{N} |X_i^2 u|^{p_i} - \sum_{i=1}^{N} X_i^2 \left(\frac{u^{p_i}}{v^{p_i-1}}\right)|X_i^2 v|^{p_i-2} X_i^2 v,
	\end{equation*}
	and
	\begin{align*}
	L_1(u,v) :=& \sum_{i=1}^{N} |X_i^2 u|^{p_i} - \sum_{i=1}^{N} p_i \left( \frac{u}{v}\right)^{p_i-1} X_i^2 u X_i^2 v |X_i^2 v|^{p_i-2} \\
	+& \sum_{i=1}^{N}(p_i-1)\left(\frac{u}{v}\right)^{p_i}|X_i^2 v|^{p_i}
	\\  -&\sum_{i=1}^{N} p_i(p_i-1) \frac{u^{p_i-2}}{v^{p_i-1}} |X_i^2 v|^{p_i-2} X_i^2 v \left(X_i u - \frac{u}{v}X_i v\right)^2.	
	\end{align*}
\end{lemma}

\begin{proof}[Proof of Lemma \ref{Picone1}]
	A direct computation gives
	\begin{align*}
	X_i^2 \left(\frac{u^{p_i}}{v^{p_i-1}}\right) &= X_i \left( p_i \frac{u^{p_i-1}}{v^{p_i-1}}X_i u - (p_i-1)\frac{u^{p_i}}{v^{p_i}}X_i v \right) \\
	& = p_i(p_i-1)\frac{u^{p_i-2}}{v^{p_i-2}} \left(\frac{(X_iu) v - u(X_i v)}{v^2}\right)X_i u + p_i \frac{u^{p_i-1}}{v^{p_i-1}} X_i^2 u \\
	& -p_i(p_i-1)\frac{u^{p_i-1}}{v^{p_i-1}}\left(\frac{(X_iu) v - u(X_i v)}{v^2}\right)X_i v - (p_i-1)\frac{u^{p_i}}{v^{p_i}}X_i^2 v\\
	&= p_i(p_i-1) \left( \frac{u^{p_i-2}}{v^{p_i-1}} |X_i u|^2 - 2\frac{u^{p_i-1}}{v^{p_i}}X_i v X_i u +\frac{u^{p_i}}{v^{p_i+1}} |X_i v|^2  \right) \\
	& + p_i \frac{u^{p_i-1}}{v^{p_i-1}} X_i^2 u - (p_i-1)\frac{u^{p_i}}{v^{p_i}} X_i^2 v \\
	& = p_i(p_i-1) \frac{u^{p_i-2}}{v^{p_i-1}} \left(X_i u - \frac{u}{v} X_i v\right)^2  + p_i \frac{u^{p_i-1}}{v^{p_i-1}} X_i^2 u - (p_i-1)\frac{u^{p_i}}{v^{p_i}} X_i^2 v,
	\end{align*}
	which gives \eqref{L=R}. By Young's inequality we have
	\begin{equation*}
	\frac{u^{p_i-1}}{v^{p_i-1}} X_i^2 u X_i^2 v |X_i^2 v|^{p_i-2} \leq \frac{|X_i^2u|^{p_i}}{p_i} + \frac{1}{q_i}\frac{u^{p_i}}{v^{p_i}}|X_i^2 v|^{p_i},\quad i=1,\ldots,N,
	\end{equation*}
	where $p_i>1, q_i>1$ $\frac{1}{p_i}+\frac{1}{q_i}=1$. Since $X_i^2 v < 0$  we arrive at
	\begin{align*}
	L_1(u,v) &\geq \sum_{i=1}^{N} |X_i^2 u|^{p_i} +\sum_{i=1}^{N}(p_i-1) \frac{u^{p_i}}{v^{p_i}} |X_i^2 v|^{p_i} - \sum_{i=1}^{N} p_i\left(\frac{|X_i^2 u|^{p_i}}{p_i}+\frac{1}{q_i} \frac{u^{p_i}}{v^{p_i}}|X_i^2 v|^{p_i}\right) \\
	& -\sum_{i=1}^{N} p_i(p_i-1)\frac{u^{p_i-2}}{v^{p_i-1}} |X_i^2 v|^{p_i-2}X^2_iv \left|X_i u - \frac{u}{v} X_iv \right|^2 \\
	&=\sum_{i=1}^{N} \left(p_i-1- \frac{p_i}{q_i}\right)\frac{u^{p_i}}{v^{p_i}}|X_i^2 v|^{p_i}
	\\& -\sum_{i=1}^{N}p_i(p_i-1) \frac{u^{p_i-2}}{v^{p_i-1}} |X_i^2 v|^{p_i-2}X_i^2v \left|X_i u -\frac{u}{v}X_iv\right|^2
	\geq 0.
	\end{align*}
	This completes the proof of Lemma \ref{Picone1}.
\end{proof}
Now we are ready to prove the subelliptic Rellich type inequality on $\mathbb{G}$.
\begin{theorem}\label{Rellich}
	Let $\Omega \subset \mathbb{G} \backslash \{x'=0\}$ be an open set, where $\mathbb{G}$ is a stratified group with $N$ being the dimension of the first stratum. Then for a function $u \geq 0 $, $u \in C^{2}(\Omega)$, and $2<\alpha_i<N-2$ we have the following inequality
	\begin{equation}\label{rellich_ineq}
	\sum_{i=1}^{N} \int_{\Omega} |X_i^2 u|^{p_i} dx \geq \sum_{i=1}^{N} C_i(\alpha_i,p_i) \int_{\Omega} \frac{|u|^{p_i}}{|x_i'|^{2p_i}}dx,
	\end{equation}
	where $1<p_i<N$ for $i=1,\ldots,N$, and
	\begin{equation*}
	C_i(\alpha_i,p_i)= (\alpha_i(\alpha_i-1))^{p_i-1} (\alpha_ip_i-2p_i-\alpha_i+2)(\alpha_ip_i-2p_i-\alpha_i+1).
	\end{equation*}
\end{theorem}

\begin{proof}[Proof of Theorem \ref{Rellich}]
	We introduce the auxiliary function
	\begin{equation*}
	v := \prod_{j=1}^{N} |x_j'|^{\alpha_j} = |x_i'|^{\alpha_i} V_i,
	\end{equation*}
	we choose $\alpha_j$ later, and let $V_i = \prod_{j=1, j\neq i}^{N}|x_j'|^{\alpha_j}$. Then we have
	\begin{align*}
	&X_i^2 v = X_i(\alpha_iV_i|x_i'|^{\alpha_i-2}x_i') =\alpha_i(\alpha_i-1)V_i |x_i'|^{\alpha_i-2},\\
	|&X_i^2 v|^{p_i-2} = (\alpha_i(\alpha_i-1))^{p_i-2}V_i^{p_i-2} |x_i'|^{\alpha_i p_i -2p_i -2\alpha_i +4}, \\
	|&X_i^2 v|^{p_i-2} X_i^2 v = (\alpha_i(\alpha_i-1))^{p_i-1}V_i^{p_i-1} |x_i'|^{\alpha_i p_i -2p_i -\alpha_i +2}.
	\end{align*}
	Consequently, we obtain 
	\begin{align*}
	X_i^2(|X_i^2 v|^{p_i-2} X_i^2 v) =& (\alpha_i(\alpha_i-1))^{p_i-1}V_i^{p_i-1} X_i^2(|x_i'|^{\alpha_i p_i -2p_i -\alpha_i +2})\\  =&(\alpha_i(\alpha_i-1))^{p_i-1}(\alpha_i p_i -2p_i -\alpha_i +2)V_i^{p_i-1}X_i\left(|x_i'|^{\alpha_i p_i -2p_i -\alpha_i}x_i'\right) \\
	=& (\alpha_i(\alpha_i-1))^{p_i-1} (\alpha_ip_i-2p_i-\alpha_i+2)(\alpha_ip_i-2p_i-\alpha_i+1) \\
	& \times V_i^{p_i-1} |x_i'|^{\alpha_i(p_i-1)-2p_i} .
	\end{align*}
	Thus, for twice differentiable function $v>0$ a.e. in $\Omega$ with $X_i^2 v < 0$ we have
	\begin{equation}\label{lem2}
	X_i^2 (|X_i^2|^{p_i-2}X_i^2 v) = C_i(\alpha_i,p_i) \frac{v^{p_i-1}}{|x'_i|^{2p_i}}
	\end{equation}
	a.e. in $\Omega$.
	Using \eqref{lem2} we compute
	\begin{align*}
	0 \leq \int_{\Omega} L_1(u,v)dx &= \int_{\Omega} R_1(u,v) dx \\
	& =	\sum_{i=1}^{N} \int_{\Omega} |X_i^2 u|^{p_i} dx - 	\sum_{i=1}^{N} \int_{\Omega} X_i^2 \left(\frac{u^{p_i}}{v^{p_i-1}}\right)|X_i^2 v|^{p_i-2}X_i^2 v dx \\
	& = \sum_{i=1}^{N} \int_{\Omega} |X_i^2 u|^{p_i} dx - 	\sum_{i=1}^{N} \int_{\Omega} \frac{u^{p_i}}{v^{p_i-1}} X_i^2\left(|X_i^2 v|^{p_i-2}X_i^2 v\right)dx\\
	& = \sum_{i=1}^{N} \int_{\Omega} |X_i^2 u|^{p_i} dx - \sum_{i=1}^{N} C_i(\alpha_i,p_i) \int_{\Omega} \frac{|u|^{p_i}}{|x_i'|^{2p_i}}dx.
	\end{align*}
	The proof of Theorem \ref{Rellich} is complete.
\end{proof}

\section{\textbf{Hardy type inequalities with multiple singularities on $\mathbb{G}$}}\label{sec3}
In this section we obtain the analogue of the Hardy inequality with multiple singularities on a stratified group. The singularities are represented by a family $\{a_k\}^m_{k=1} \in \mathbb{G},$ where we write $a_k= (a_k',a_k''),$ with $a_k'$ being in the first stratum of $\mathbb{G}$. We can also write $a_k'=(a_{k1}',\ldots,a_{kN}')$, so $(xa^{-1}_k)'=x'-a_k'.$ 
\begin{theorem}\label{thm_L}
	Let $\Omega \subset \mathbb{G}$ be an open set, where $\mathbb{G}$ is a stratified group  with $N$ being the dimension of the first stratum. Let $N \geq 3$, $x = (x',x'') \in \mathbb{G}$ with $x'=(x'_1,\ldots,x'_N)$ being in the first stratum of $\mathbb{G}$, and let $ a_k\in \mathbb{G},  k=1,\ldots,m,$ be the singularities. Then we have
	\begin{equation}\label{L2.4}
	\int_{\Omega} |\nabla_{\mathbb{G}} u |^2 dx \geq\left(\frac{N-2}{2}\right)^2 \int_{\Omega} \frac{\sum_{j=1}^{N}\left|\sum_{k=1}^{m}\frac{(x a_{k}^{-1})_j'}{|(xa_k^{-1})'|^N}\right|^2}{\left(\sum_{k=1}^{m} \frac{1}{|(xa_k^{-1})'|^{N-2}}\right)^2} |u|^2 dx,
	\end{equation}
	for all $u \in C_0^{\infty}(\Omega)$.
\end{theorem}
The Euclidean case of this inequality was obtained by Kapitanski and Laptev  \cite{KL}. In \eqref{L2.4}, $(x a_k^{-1})'_j= x_j'-a_{kj}'$ denotes the $j^{\text{th}}$ component of $xa_k^{-1}$.

\begin{proof}[Proof of Theorem \ref{thm_L}]
	Let us introduce a vector-field $\mathcal{A}(x)=(\mathcal{A}_{1}(x),\ldots, \mathcal{A}_{N}(x))$ to be specified later. Also let $\lambda$ be a real parameter for optimisation. We start with the inequality
	\begin{align*}
	0 \leq& \int_{\Omega} \sum_{j=1}^{N} (|X_j u-\lambda \mathcal{A}_j u|^2)dx  \\
	=& \int_{\Omega} \left( |\nabla_{\mathbb{G}} u|^2 - 2\lambda \text{Re} \sum_{j=1}^{N} \overline{\mathcal{A}_j u} X_j u + \lambda^2  \sum_{j=1}^{N} |\mathcal{A}_j|^2 |u|^2 \right)dx.
	\end{align*}
	By using the integration by parts we get 
	\begin{equation}\label{sup}
	- \int_{\Omega} \left(\lambda^2  \sum_{j=1}^{N} |\mathcal{A}_j|^2 + \lambda {\rm div}_{\mathbb{G}} \mathcal{A} \right)|u|^2dx \leq \int_{\Omega} |\nabla_{\mathbb{G}} u|^2 dx.
	\end{equation}
	We differentiate the integral on the left-hand side with respect to $\lambda$ to optimise it, yielding 
	\begin{equation*}
	2\lambda |\mathcal{A}|^2 + {\rm div}_{\mathbb{G}} {\mathcal{A}} = 0,
	\end{equation*}
	for all $x \in \Omega$. This is a restriction on $\mathcal{A}(x)$ giving $ \frac{{\rm div}_{\mathbb{G}} \mathcal{A}(x)}{|\mathcal{A}(x)|^2} = const.$ For $\lambda = \frac{1}{2}$ we get
	\begin{equation}\label{L2.1}
	{\rm div}_{\mathbb{G}} \mathcal{A}(x) = - |\mathcal{A}(x)|^2.
	\end{equation}
	Then putting \eqref{L2.1} in \eqref{sup} we have the following Hardy inequality
	\begin{equation*}
	\frac{1}{4} \int_{\Omega} \sum_{j=1}^{N} |\mathcal{A}_j(x)|^2 |u|^2 dx \leq \int_{\Omega} |\nabla_{\mathbb{G}} u|^2 dx.
	\end{equation*}
	Now if we assume that $\mathcal{A} = \nabla_{\mathbb{G}} \phi $ for some function $\phi$, then \eqref{L2.1} becomes
	\begin{equation*}
	\mathcal{L} \phi + |\nabla_{\mathbb{G}} \phi |^2 = 0.
	\end{equation*}
	It follows that the function is harmonic (with respect to the sub-Laplacian $\mathcal{L}$).
	\begin{equation*}
	w = e^{\phi} \geq 0
	\end{equation*}
	Then $w$ is a constant $>0$ or has a singularity. Let us consider
	\begin{equation*}
	w := \sum_{k=1}^{m} \frac{1}{|(xa_k^{-1})'|^{N-2}},
	\end{equation*}
	and then take
	\begin{equation*}
	\phi (x) = \ln (w).
	\end{equation*}
	Therefore
	\begin{align*}
	\mathcal{A}(x) =  \nabla_{\mathbb{G}}(\ln w) =& \frac{1}{w} \nabla_{\mathbb{G}} \left(\sum_{k=1}^{m} |(xa_k^{-1})'|^{2-N}\right) \\
	=& \frac{1}{w} \sum_{k=1}^{m} \nabla_{\mathbb{G}} \left(\sum_{j=1}^{N} ((x a_{k}^{-1})_j')^2\right)^{\frac{2-N}{2}}
	\\=& -\frac{N-2}{w} \left( \sum_{k=1}^{m} \frac{(xa_k^{-1})'}{|(xa_k^{-1})'|^N} \right),
	\end{align*}
	and
	\begin{equation*}
	|\mathcal{A}(x)|^2 = \sum_{j=1}^{N} |\mathcal{A}_j (x)|^2 = \left(\frac{N-2}{w}\right)^2 \sum_{j=1}^{N}\left|\sum_{k=1}^{m}\frac{(xa_{k}^{-1})_j'}{|(xa_k^{-1})'|^N}\right|^2.
	\end{equation*}
	This completes the proof of Theorem \ref{thm_L}.
\end{proof}
We then also obtain the corresponding uncertainty principle.
\begin{corollary}\label{uncert1}
	Let $\Omega \subset \mathbb{G}$ be an open set, where $\mathbb{G}$ is a stratified group  with $N$ being the dimension of the first stratum. Let $N \geq 3$, $x = (x',x'') \in \mathbb{G}$ with $x'=(x'_1,\ldots,x'_N)$ being in the first stratum of $\mathbb{G}$. Let $ a_k\in \mathbb{G}, k=1,\ldots,m,$ be the singularities. Then we have
	\begin{equation}
	\frac{N-2}{2} \int_{\Omega} |u|^2 dx \leq \left(\int_{\Omega} |\nabla_{\mathbb{G}} u|^2 dx\right)^{\frac{1}{2}} \left(\int_{\Omega} \frac{\left(\sum_{k=1}^{m} \frac{1}{|(xa_k^{-1})'|^{N-2}}\right)^2}{\sum_{j=1}^{N}\left|\sum_{k=1}^{m}\frac{(xa_{k}^{-1})_j' }{|(xa_k^{-1})'|^N}\right|^2} |u|^2 dx\right)^{\frac{1}{2}},
	\end{equation}
	for all $u \in C_0^{\infty}(\Omega)$ and $1<p_i<N$ for $i=1,\ldots,N.$
\end{corollary}
\begin{proof}[Proof of Corollary \ref{uncert1}]  By \eqref{L2.4} and the Cauchy-Schwarz inequality we get
	\begin{align*}
	&\int_{\Omega} |\nabla_{\mathbb{G}} u|^2 dx \int_{\Omega} \frac{\left(\sum_{k=1}^{m} \frac{1}{|(xa_k^{-1})'|^{N-2}}\right)^2}{\sum_{j=1}^{N}\left|\sum_{k=1}^{m}\frac{(xa_{k}^{-1})_j'}{|(xa_k^{-1})'|^N}\right|^2} |u|^2 dx \\
	& \geq \left(\frac{N-2}{2}\right)^2 \int_{\Omega} \frac{\sum_{j=1}^{N}\left|\sum_{k=1}^{m}\frac{ (xa_{k}^{-1})_j'}{|(xa_k^{-1})'|^N}\right|^2}{\left(\sum_{k=1}^{m} \frac{1}{|(xa_k^{-1})'|^{N-2}}\right)^2} |u|^2 dx \int_{\Omega} \frac{\left(\sum_{k=1}^{m} \frac{1}{|(xa_k^{-1})'|^{N-2}}\right)^2}{\sum_{j=1}^{N}\left|\sum_{k=1}^{m}\frac{(xa_{k}^{-1})_j'}{|(xa_k^{-1})'|^N}\right|^2} |u|^2 dx \\
	& \geq \left(\frac{N-2}{2} \right)^2 \left( \int_{\Omega}|u|^2 dx\right)^{2}.
	\end{align*}
	The proof is complete.
\end{proof}

\section{\textbf{Many-particle Hardy type inequality on $\mathbb{G}$}}\label{sec4}

Suppose now that we have $n$ particles, where $n$ is a positive integer. Let $\mathbb{G}^n$ be the product $\mathbb{G}^n := \overbrace{\mathbb{G}\times \ldots \times \mathbb{G}}^n$. Now let us consider $x=(x_1,\ldots,x_n) \in \mathbb{G}^n,$ $ x_j \in \mathbb{G} $. Let $x \in \mathbb{G}^n$ with  $x'=(x'_1,\ldots,x'_n)$ and $x'_i=(x'_{i1},\ldots,x'_{iN})$ being the coordinates on the first stratum of $\mathbb{G}$ for $i=1,\ldots,n$. The distance between particles $x_i,x_j \in \mathbb{G}$ can be defined by
\begin{equation*}
r_{ij} := |(x_ix_j^{-1})'|=|x'_i-x'_j| = \sqrt{\sum_{k=1}^{N} (x'_{ik}-x'_{jk})^2}.
\end{equation*}
We also use the following notation
\begin{equation*}
\nabla_{\mathbb{G}_i} = (X_{i1},\ldots,X_{iN})
\end{equation*}
for the gradient associated to the $i$-th particle. We denote $\nabla_{\mathbb{G}^n}:= (\nabla_{\mathbb{G}_1},\ldots, \nabla_{\mathbb{G}_n})$, and
\begin{equation*}
\mathcal{L}_i = \sum_{k=1}^{N}X_{ik}^2
\end{equation*}
is the sub-Laplacian associated to the $i$-th particle. We note that $\mathcal{L} = \sum_{i=1}^{N} \mathcal{L}_i.$
We are new ready to prove the following crucial inequality in $\mathbb{R}^m.$
\begin{lemma}\label{lem_3.1}
	Let $m\geq 1$, and let
	\begin{equation*}
	\mathcal{A} = \left(\mathcal{A}_1(x),\ldots,\mathcal{A}_m(x)\right)
	\end{equation*}
	be a mapping in $\mathcal{A}:\mathbb{R}^m \rightarrow \mathbb{R}^m$ whose components and their first derivatives are uniformly bounded in $\mathbb{R}^m$. Then for $u \in C_0^1(\mathbb{R}^m)$ we have
	\begin{equation}\label{lem_eq}
	\int_{\mathbb{R}^m} |\nabla u|^2 dx \geq \frac{1}{4} \frac{\left(\int_{\mathbb{R}^m} {\rm div} \mathcal{A} |u|^2 dx\right)^2}{\int_{\mathbb{R}^m}|\mathcal{A}|^2|u|^2dx}.
	\end{equation}
\end{lemma}
\begin{proof}[Proof of Lemma \ref{lem_3.1}]
	We have
	\begin{align*}
	\left|\int_{\mathbb{R}^m} {\rm div} \mathcal{A} |u|^2dx\right| &= 2\left|\text{Re} \int_{\mathbb{R}^m} \langle \mathcal{A},\nabla u\rangle\overline{u}dx\right| \\
	& \leq 2 \left(\int_{\mathbb{R}^m}|\mathcal{A}|^2 |u|^2dx\right)^{1/2} \left(\int_{\mathbb{R}^m}|\nabla u|^2dx\right)^{1/2}.
	\end{align*}
	We have used the Cauchy-Schwarz inequality in the last line. The proof is finished by squaring this inequality.
\end{proof}

\begin{theorem}\label{thm1}
	Let $\Omega \subset \mathbb{G}^n$ be an open set, where $\mathbb{G}$ is a stratified group  with $N$ being the dimension of the first stratum. Let $N\geq 2$ and $n\geq 3$. Let $r_{ij} =|(x_ix_j^{-1})'|= |x'_i-x'_j|.$ Then we have
	
	\begin{align}\label{H1}
	\int_{\Omega} |\nabla_{\mathbb{G}^n} u|^2 dx \geq \frac{(N-2)^2}{n} \int_{\Omega}\sum_{1 \leq i<j\leq n}\frac{|u|^2}{r_{ij}^2} dx,
	\end{align}
	for all $u \in C^1(\Omega)$.
\end{theorem}

The Euclidean case of inequality  \eqref{H1} is obtained by M. Hoffmann-Ostenhof, T. Hoffmann-Ostenhof, A. Laptev, and J. Tidblom  \cite{HHLT}.
\begin{proof}[Proof of Theorem \ref{thm1}]
	Let us choose a mapping $ \mathcal{B}_1$ in the following form
	\begin{equation*}
	\mathcal{B}_1(x'_i,x'_j) := \frac{(x_ix_j^{-1})'}{r^2_{ij}}, \quad 1\leq i < j\leq n.
	\end{equation*}
	And putting the mapping $\mathcal{B}_1$ in \eqref{lem_eq} we have
	\begin{align}\label{4.1}
	\int_{\Omega} |(\nabla_{\mathbb{G}_i} -\nabla_{\mathbb{G}_j})u|^2 dx &\geq \frac{1}{4} \frac{\left(\int_{\Omega}\left(({\rm div}_{\mathbb{G}_i}-{\rm div}_{\mathbb{G}_j})\mathcal{B}_1\right)|u|^2dx\right)^2}{\int_{\Omega}|\mathcal{B}_1|^2|u|^2dx} \nonumber \\
	& = \frac{1}{4} \frac{\left(\int_{\Omega} \frac{2(N-2)}{|(x_ix_j^{-1})'|^2}|u|^2dx\right)^2}{\int_{\Omega} \frac{|u|^2}{|(x_ix_j^{-1})'|^2}dx}
	\nonumber \\
	&= (N-2)^2 \int_{\Omega} \frac{|u|^2}{r^2_{ij}} dx.
	\end{align}
	Also, we introduce another mapping $\mathcal{B}_2$
	\begin{equation*}
	\mathcal{B}_2(x) := \frac{\sum_{j=1}^{n}x'_j}{\left|\sum_{j=1}^{n}x'_j\right|^2},
	\end{equation*}
	and 
	\begin{align*}
	{\rm div_{\mathbb{G}_i}} \mathcal{B}_2 = \nabla_{\mathbb{G}_i} \cdot \mathcal{B}_2 &= \sum_{k=1}^{N} X_{ik} \left( \frac{\sum_{j=1}^{n} x_{jk}'}{|\sum_{j=1}^{n} x_j'|^2}  \right) \\
	& = \frac{Nn |\sum_{j=1}^{n} x_j'|^2 - 2n\left((\sum_{j=1}^{n} x_{j1}')^2 + \ldots+ (\sum_{j=1}^{n} x_{jN}')^2 \right)}{|\sum_{j=1}^{n} x_j'|^4} \\
	& = \frac{Nn-2n}{|\sum_{j=1}^{n} x_j'|^2}.
	\end{align*}
	As before we put the mapping $\mathcal{B}_2$ in \eqref{lem_eq} and using above computation yielding 
	\begin{align}\label{4.2}
	\int_{\Omega} \left|\sum_{i=1}^{n} \nabla_{\mathbb{G}_i} u \right|^2 dx &\geq \frac{1}{4} \frac{\left(\int_{\Omega} (\sum_{i=1}^{n}{\rm div }_{\mathbb{G}_i}\mathcal{B}_2)|u|^2dx\right)^2}{\int_{\Omega}|\mathcal{B}_2|^2|u|^2dx} \nonumber \\ 
	& = \frac{1}{4} \frac{\left(\int_{\Omega} \sum_{i=1}^{n}\frac{Nn-2n}{|\sum_{j=1}^{n}x_j'|^2}|u|^2dx\right)^2}{\int_{\Omega} \frac{|u|^2}{\left|\sum_{j=1}^{n}x'_j\right|^2}dx}\nonumber \\
	&=\frac{(N-2)^2n^4}{4} \int_{\Omega} \frac{|u|^2}{\left|\sum_{j=1}^{n}x'_j\right|^2}dx.
	\end{align}
	Adding inequalities \eqref{4.1} and \eqref{4.2} and using the identity
	\begin{equation*}
	n \sum_{i=1}^n |\nabla_{\mathbb{G}_i} u|^2 = \sum_{1\leq i<j\leq n} \left|\nabla_{\mathbb{G}_i}u - \nabla_{\mathbb{G}_j}u\right|^2 + \left|\sum_{i=1}^n \nabla_{\mathbb{G}_i}u\right|^2,
	\end{equation*}
	we arrive at
	\begin{align}\label{4.3}
	\sum_{i=1}^n\int_{\Omega} |\nabla_{\mathbb{G}_i} u|^2 dx \geq \frac{(N-2)^2}{n} \int_{\Omega} \sum_{i<j} \frac{|u|^2}{r^2_{ij}} dx + \frac{(N-2)^2n^3}{4}\int_{\Omega} \frac{|u|^2}{\left|\sum_{j=1}^{n}x'_j\right|^2}dx.
	\end{align}
	Because the last term on right-hand side is positive, we get
	\begin{align*}
	\sum_{i=1}^n \int_{\Omega} |\nabla_{\mathbb{G}_i} u|^2 dx \geq \frac{(N-2)^2}{n} \int_{\Omega}\sum_{i<j}\frac{|u|^2}{r_{ij}^2} dx.
	\end{align*}
	Also we have 
	\begin{align*}
	\sum_{i=1}^{n} |\nabla_{\mathbb{G}_i} u|^2 
	& = (\nabla_{\mathbb{G}_1} u)^2+\ldots+ (\nabla_{\mathbb{G}_n}u)^2\\
	&= |(\nabla_{\mathbb{G}_1} u,\ldots, \nabla_{\mathbb{G}_n}u)|^2 \\
	& = 	|\nabla_{\mathbb{G}^n} u |^2.		
	\end{align*}
	The proof of Theorem \ref{thm1} is complete.
\end{proof}

The following theorem deals with the total separation of $n \geq 2$ particles.
\begin{theorem}\label{thm22} Let $\Omega \subset \mathbb{G}^n$ be an open set, where $\mathbb{G}$ is a stratified group  with $N$ being the dimension of the first stratum. Let $\rho^2 := \sum_{i<j}| (x_ix_j^{-1})'|^2 =\sum_{i<j} |x_i'-x_j'|^2 $ with $x'_i \neq x'_j$. Then we have
	\begin{equation}\label{Hardy_d}
	\int_{\Omega} |\nabla_{\mathbb{G}} u|^2 dx = n\left(\frac{(n-1)}{2}N-1\right)^2 \int_{\Omega} \frac{|u|^2}{\rho^2} dx + \int_{\Omega} |\nabla_{\mathbb{G}}\rho^{-2\alpha}u|^2 \rho^{4\alpha} dx
	\end{equation}
	for all $u \in C_0^{\infty}(\Omega)$ with $\alpha = \frac{2-(n-1)N}{4}$.
\end{theorem}
The Euclidean case of inequality \eqref{Hardy_d} was obtained by Douglas Lundholm \cite{Lundholm}.
\begin{proposition}\label{prop1}
	Let $\Omega \subset \mathbb{G}^n$ be an open set, where $\mathbb{G}$ is a stratified group  with $N$ being the dimension of the first stratum. Let $f:\Omega \rightarrow  (0,\infty)$ be twice differentiable. Then for any function $u \in C_0^{\infty}(\Omega)$ and $\alpha \in \mathbb{R}$, we have
	\begin{equation}\label{eq3}
	\int_{\Omega} |\nabla_{\mathbb{G}}u|^2 dx = \int_{\Omega} \left( \alpha(1-\alpha)\frac{|\nabla_{\mathbb{G}} f|^2}{f^2} - \alpha \frac{ \mathcal{L}f}{f}\right)|u|^2 dx + \int_{\Omega} |\nabla_{\mathbb{G}} v|^2 f^{2\alpha} dx,
	\end{equation}
	where $v := f^{-\alpha}u$.
\end{proposition}
\begin{proof}[Proof of Proposition \ref{prop1}]
	Let us compute for $u = f^{\alpha}v$, that 
	$$\nabla_{\mathbb{G}} u = \alpha f^{\alpha -1}(\nabla_{\mathbb{G}} f) v + f^{\alpha}\nabla_{\mathbb{G}} v.$$
	Then by squaring the above expression we have
	\begin{align*}
	|\nabla_{\mathbb{G}} u|^2 &= \alpha^2 f^{2(\alpha-1)} |\nabla_{\mathbb{G}} f|^2 |v|^2 + \text{Re}(2 \alpha v f^{2\alpha-1}(\nabla_{\mathbb{G}} f) \cdot  (\nabla_{\mathbb{G}} v))  + f^{2\alpha}|\nabla_{\mathbb{G}} v|^2 \\
	& =\alpha^2 f^{2(\alpha-1)} |\nabla_{\mathbb{G}} f|^2 |v|^2 +\alpha f^{2\alpha-1}(\nabla_{\mathbb{G}} f) \cdot  \nabla_{\mathbb{G}} |v|^2 + f^{2\alpha}|\nabla_{\mathbb{G}} v|^2.
	\end{align*}
	By integrating this expression over $\Omega$, we have 
	\begin{align*}
	\int_{\Omega}|\nabla_{\mathbb{G}} u|^2 dx &=  \int_{\Omega} \alpha^2 f^{2(\alpha-1)} |\nabla_{\mathbb{G}} f|^2 |v|^2dx \\&+ \int_{\Omega}\text{Re}(\alpha f^{2\alpha-1}(\nabla_{\mathbb{G}} f) \cdot  \nabla_{\mathbb{G}} |v|^2) dx
	+ \int_{\Omega} f^{2\alpha}|\nabla_{\mathbb{G}} v|^2dx \\
	&= \int_{\Omega} \alpha^2 f^{2(\alpha-1)} |\nabla_{\mathbb{G}} f|^2 |v|^2dx \\&
	-\alpha \int_{\Omega} \nabla_{\mathbb{G}} \cdot (f^{2\alpha-1}\nabla_{\mathbb{G}} f) |v|^2 dx + \int_{\Omega} f^{2\alpha}|\nabla_{\mathbb{G}} v|^2dx.
	\end{align*}
	We have used the integration by parts to the middle term on the right-hand side. Then
	\begin{equation*}
	\nabla_{\mathbb{G}} \cdot (f^{2 \alpha -1}\nabla_{\mathbb{G}} f) = (2\alpha -1)f^{2\alpha-2} |\nabla_{\mathbb{G}} f|^2 + f^{2\alpha-1} \mathcal{L} f,
	\end{equation*}
	and by using this fact we get
	\begin{align*}
	\int_{\Omega}|\nabla_{\mathbb{G}} u|^2 dx& = \int_{\Omega} \alpha^2 f^{2(\alpha-1)} |\nabla_{\mathbb{G}} f|^2 |v|^2dx - \int_{\Omega} \alpha f^{2\alpha-1} \mathcal{L} f |v|^2 dx \\
	& - \int_{\Omega} \alpha (2\alpha -1)f^{2\alpha-2} |\nabla_{\mathbb{G}} f|^2 |v|^2 dx + \int_{\Omega} f^{2\alpha}|\nabla_{\mathbb{G}} v|^2dx.
	\end{align*}
	Putting back $v = f^{-\alpha} u$ and collecting the terms we arrive at \eqref{eq3}.
\end{proof}
\begin{proof}[Proof of Theorem \ref{thm22}]
	The following computation gives
	\begin{equation*}
	\nabla_{\mathbb{G}_k} \rho^2 = (X_{k1} \rho^2,\ldots,X_{kN} \rho^2) =2 \sum_{k\neq j}^n (x_kx_j^{-1})',
	\end{equation*}
	where $\nabla_{\mathbb{G}_k} = (X_{k1},\ldots,X_{kN})$.
	Hence
	\begin{align}\label{Lp}
	\mathcal{L} \rho^2 
	=2 \sum_{k=1}^n \sum_{k\neq j}^n \nabla_{\mathbb{G}_k} \cdot (x_kx_j^{-1})' = 2n(n-1)N,
	\end{align}
	\begin{equation}\label{r^2}
	|\nabla_{\mathbb{G}} \rho^2|^2 = 8 \sum_{1\leq i<j\leq n} |(x_kx_j^{-1})'|^2 + 8 \sum_{k=1}^n \sum_{1\leq i<j\leq n} (x_kx_i^{-1})'\cdot(x_kx_j^{-1})' = 4n\rho^2,
	\end{equation}
	where in the last step we used the identity
	\begin{equation*}
	\sum_{k=1}^n \sum_{1\leq i<j\leq n} (x_k x_i^{-1})' \cdot (x_k x_j^{-1})' = \frac{n-2}{2} \sum_{1\leq i<j\leq n} |(x_i x_j^{-1})'|^2.
	\end{equation*}
	By putting \eqref{Lp} and \eqref{r^2} in Proposition \ref{prop1} with $f = \rho^2$ we have
	\begin{equation*}
	\int_{\Omega} |\nabla_{\mathbb{G}}u|^2 dx = 4n\alpha \left( \frac{2-(n-1)N}{2} - \alpha\right)\int_{\Omega}\frac{|u|^2}{\rho^2}dx +   \int_{\Omega} |\nabla_{\mathbb{G}} \rho^{-2\alpha}u|^2 \rho^{4\alpha} dx.
	\end{equation*}
	To optimise we differentiate the integral 
	\begin{equation*}
	4n\alpha \left( \frac{2-(n-1)N}{2} - \alpha\right)\int_{\Omega}\frac{|u|^2}{\rho^2}dx
	\end{equation*}
	with respect to $\alpha$, then we have
	\begin{equation*}
	\frac{2-(n-1)N}{2} -2\alpha =0,
	\end{equation*} 
	and
	\begin{equation*}
	\alpha = \frac{2-(n-1)N}{4},
	\end{equation*}
	which completes the proof of Theorem \ref{thm22}.
\end{proof}

\section{\textbf{Hardy type inequalities with exponential weights on $\mathbb{G}$}}\label{sec5}

In this section, we get the horizontal Hardy inequality with exponential weights on $\mathbb{G}$.

\begin{theorem}\label{lemma2}
	Let $\Omega \subset \mathbb{G} $ be an open set, where $\mathbb{G}$ is a stratified group with $ N\geq 3$ being the dimension of the first stratum.
	Let $x_0 \in \Omega$. Then  we have
	
	\begin{align}\label{PH2.9}
	\int_{\Omega} e^{-\frac{|(xx_0^{-1})'|^2}{4\lambda}} \left( \frac{(N-2)^2}{4|x'|^2} - \frac{N}{4\alpha} + \frac{|(xx_0^{-1})'|^2}{16\lambda^2} \right) |u|^2 dx \leq \int_{\Omega} e^{-\frac{|(xx_0^{-1})'|^2}{4\lambda}}  |\nabla_{\mathbb{G}} u|^2 dx
	\end{align}
	for all $u \in C^1(\Omega)$ and for each $ \lambda>0$.
\end{theorem}
In the Euclidean case, this inequality is called two parabolic-type Hardy inequality, which was obtained by Zhang \cite{Zhang}.
\begin{proof}[Proof of Theorem \ref{lemma2}]
	Let us recall the horizontal Hardy inequality \cite{RS17a} for all $v \in C^1 (\Omega),$ 
	\begin{equation}\label{PH2.10}
	\frac{(N-2)^2}{4} \int_{\Omega} \frac{|v|^2}{|x'|^2} dx \leq \int_{\Omega} |\nabla_{\mathbb{G}} v|^2 dx.
	\end{equation}
	Let $v =  e^{-\frac{|(xx_0^{-1})'|^2}{8\lambda}} u$, then
	\begin{equation*}
	\nabla_{\mathbb{G}} v =  e^{-\frac{|(xx_0^{-1})'|^2}{8\lambda}} \nabla_{\mathbb{G}}u - \frac{(xx_0^{-1})'}{4\lambda}  e^{-\frac{|(xx_0^{-1})'|^2}{8\lambda}} u,
	\end{equation*}
	for all $v \in C^1(\Omega).$ Then by inequality \eqref{PH2.10} we have
	\begin{align}\label{6.3}
	\frac{(N-2)^2}{4} \int_{\Omega}  e^{-\frac{|(xx_0^{-1})'|^2}{4\lambda}} \frac{|u|^2}{|x'|^2}dx \leq &\int_{\Omega}  e^{-\frac{|(xx_0^{-1})'|^2}{4\lambda}}|\nabla_{\mathbb{G}}u|^2 + \frac{|(xx_0^{-1})'|^2}{16 \lambda^2}  e^{-\frac{|(xx_0^{-1})'|^2}{4\lambda}}|u|^2dx \\
	& - \text{Re}\frac{1}{2\lambda} \int_{\Omega} (xx_0^{-1})' \cdot (\nabla_{\mathbb{G}}u) u e^{-\frac{|(xx_0^{-1})'|^2}{4\lambda}}dx. \nonumber
	\end{align}
	By the integration by parts in the last term of right-hand side of the inequality we have
	\begin{equation}\label{PH2.11}
	\text{Re}	\int_{\Omega}(xx_0^{-1})' \cdot (\nabla_{\mathbb{G}}u) u e^{-\frac{|(xx_0^{-1})'|^2}{4\lambda}} dx = - \frac{1}{2} \int_{\Omega} \left(N- \frac{|(xx_0^{-1})'|^2}{2\lambda} \right)  e^{-\frac{|(xx_0^{-1})'|^2}{4\lambda}} |u|^2dx.
	\end{equation}
	By putting equality \eqref{PH2.11} in \eqref{6.3} and rearranging it, we prove Theorem \ref{lemma2}.
\end{proof}

{\bf Acknowledgments.} The first
author was supported by the EPSRC Grant 
EP/R003025/1, by the Leverhulme Research Grant RPG-2017-151, and by the FWO Odysseus grant.  The second author was supported by the MESRK target program BR05236656 and the Nazarbayev University SPG grant SST 2018040. The third author was supported in parts by the Nazarbayev University SPG grant. No new data was collected or generated during the course of this research.

\bibliographystyle{amsplain}

\end{document}